DISTRIBUTIONS OF ROOTS OF REDUCED
CUBIC EQUATIONS WITH RANDOM COEFFICIENTS

A Thesis

Submitted to the Graduate Faculty of the
Louisiana State University and
Agricultural and Mechanical College
in partial fulfillment of the
requirements for the degree of
Master of Science

in

The Department of Mathematics

by
Kerry Michael Soileau
B.S., Florida Technological University, 1976
May 1980

# ACKNOWLEDGEMENT

The author wishes to express his sincere appreciation to Professor Hui-Hsiung Kuo for his help and guidance in the preparation of this thesis.

# TABLE OF CONTENTS

ABSTRACT

If the coefficients of polynomials are selected by some random process, the zeros of the resulting polynomials are in some sense random. In this paper the author rephrases the above in more precise language, and calculates the joint conditional densities of a random vector whose values determine almost surely the zeros of a "random" reduced cubic.

INTRODUCTION

<u>General Problem</u>: We consider the random function $\varphi:\mathbb{C}\times\Omega\to\mathbb{C}$ given by $\varphi(z,\omega)=z^3+\xi_1(\omega)z+\xi_2(\omega)$, where $(\Omega,\Im,P)$ is a probability space, and $\xi_1,\xi_2:\Omega\to\mathbb{R}$ are real random variables on $\Omega$ with continuous joint density $f:\mathbb{R}\times\mathbb{R}\to\mathbb{R}$. We can view $\varphi$ as a function which associates with each point $\omega\in\Omega$ a reduced cubic with real coefficients. A well-known result in the theory of equations gives the following information about the roots of $z^3+\xi_1(\omega)z+\xi_2(\omega)$:

if $\dfrac{\xi_2^2(\omega)}{4}+\dfrac{\xi_1^3(\omega)}{27}>0$, then there are exactly one real root and two conjugate imaginary roots; if $\dfrac{\xi_2^2(\omega)}{4}+\dfrac{\xi_1^3(\omega)}{27}=0$, then there are three real roots, two or more of which are equal. If $\dfrac{\xi_2^2(\omega)}{4}+\dfrac{\xi_1^3(\omega)}{27}<0$, then there are three distinct real roots. These three conditions define a disjoint partition of $\Omega$ into three events, $\mathcal{D}$, $\mathcal{S}$, and $\mathcal{K}$ respectively. Since $\mathcal{S}$ is a zero-probability event, we will neglect it and be concerned with only with the events $\mathcal{D}$ and $\mathcal{K}$, both of which we assume have nonzero probability. We will then calculate the densities $h(x,y|\mathcal{D})$ and $h(x,y|\mathcal{K})$, corresponding to the joint conditional densities of $R^*$ relative to the hypotheses $\mathcal{D}$ and $\mathcal{K}$, respectively, where $R^*:\Omega\to\mathbb{R}\times\mathbb{R}$ is a random vector on $\Omega$ whose values for each $\omega\in\Omega$ determine almost surely the zeros of the polynomial $z^3+\xi_1(\omega)z+\xi_2(\omega)$.

<u>Relation to Existing Literature</u>: This thesis is an extension of a paper by John W. Hamblen [5] which describes the solution of a similar problem in which the random function $\varphi:\mathbb{C}\times\Omega\to\mathbb{C}$ is given by $\varphi(z,\omega)=z^2+\xi_1(\omega)z+\xi_2(\omega)$.

<u>Definition of the Random Vector $R^*$</u>

We define $R^*:\Omega\to\mathbb{R}\times\mathbb{R}$

$$R^*:\omega\mapsto\left(R_1^*(\omega),R_2^*(\omega)\right)$$

to be the function given by

$$R_1^*(\omega)=\begin{cases}\{\mathrm{Re}(z);\varphi(z,\omega)=0, z\notin\mathbb{R}\} & \text{if } \dfrac{\xi_2^2(\omega)}{4}+\dfrac{\xi_1^3(\omega)}{27}>0 \\[2mm] \max\{z;\varphi(z,\omega)=0\} & \text{if } \dfrac{\xi_2^2(\omega)}{4}+\dfrac{\xi_1^3(\omega)}{27}<0\end{cases}$$

$$R_2^*(\omega)=\begin{cases}\{|\mathrm{Im}(z)|;\varphi(z,\omega)=0; z\notin\mathbb{R}\} & \text{if } \dfrac{\xi_2^2(\omega)}{4}+\dfrac{\xi_1^3(\omega)}{27}>0 \\[2mm] \max\{z;\varphi(z,\omega)=0, z<R_1^*(\omega)\} & \text{if } \dfrac{\xi_2^2(\omega)}{4}+\dfrac{\xi_1^3(\omega)}{27}<0\end{cases}$$

In short, if we choose $\omega\in\Omega$ and $\omega$ is such that $z^3+\xi_1(\omega)z+\xi_2(\omega)=0$ has one real and two complex conjugate zeros, then we define the first and second coordinates of $R^*(\omega)$ to be the real part and the absolute value of the imaginary parts of the complex conjugate zeros. If for our choice of $\omega$ the polynomial $z^3+\xi_1(\omega)z+\xi_2(\omega)=0$ has three distinct real zeros, then we define the first and second coordinates of $R^*(\omega)$ to be the largest and the next largest of these zeros. These two requirements define the value of $R^*(\omega)$ for all $\omega \in \Omega$ such that $\dfrac{\xi_2^2(\omega)}{4}+\dfrac{\xi_1^3(\omega)}{27} \neq 0$, so $R^*$ is defined almost surely on $\Omega$.

# CALCULATION OF CONDITIONAL DENSITY $h(x,y|\mathcal{D})$

The calculation of the conditional density $h(x,y|\mathcal{D})$ can be made much simpler by defining random variables $A$ and $B$ with respect to the event $\mathcal{D}$. Then the functional relation between the random variables $A$ and $B$ and the random variables $R_1^*$ and $R_2^*$ is exploited to calculate the joint conditional density of $R_1^*$ and $R_2^*$ with respect to $\mathcal{D}$.

The Joint Density $g_{A,B}(x,y|\mathcal{D})$.

If $\omega \in \Omega$ is such that $\dfrac{\xi_2^2(\omega)}{4} + \dfrac{\xi_1^3(\omega)}{27} > 0$, then the general solution of the cubic equation yields

$$R_1^*(\omega) = -\frac{1}{2}\left\{ \sqrt[3]{-\frac{\xi_2(\omega)}{2} + \sqrt{\frac{\xi_2^2(\omega)}{4} + \frac{\xi_1^3(\omega)}{27}}} + \sqrt[3]{-\frac{\xi_2(\omega)}{2} - \sqrt{\frac{\xi_2^2(\omega)}{4} + \frac{\xi_1^3(\omega)}{27}}} \right\}$$

and

$$R_2^*(\omega) = \frac{\sqrt{3}}{2}\left\{ \sqrt[3]{-\frac{\xi_2(\omega)}{2} + \sqrt{\frac{\xi_2^2(\omega)}{4} + \frac{\xi_1^3(\omega)}{27}}} - \sqrt[3]{-\frac{\xi_2(\omega)}{2} - \sqrt{\frac{\xi_2^2(\omega)}{4} + \frac{\xi_1^3(\omega)}{27}}} \right\}$$

Let $A(\omega) = \dfrac{1}{\sqrt{3}} R_2^*(\omega) - R_1^*(\omega)$ and $B(\omega) = -\dfrac{1}{\sqrt{3}} R_2^*(\omega) - R_1^*(\omega)$. We may calculate the joint density of $A$ and $B$ under the hypothesis $\dfrac{\xi_2^2(\omega)}{4} + \dfrac{\xi_1^3(\omega)}{27} > 0$, which was denoted by $\mathcal{D}$. In the following we assume $0 < P(\mathcal{D}) < 1$. We may write $A$ and $B$ explicitly as

$$A(\omega) = \sqrt[3]{-\frac{\xi_2(\omega)}{2} + \sqrt{\frac{\xi_2^2(\omega)}{4} + \frac{\xi_1^3(\omega)}{27}}} \quad \text{and} \quad B(\omega) = \sqrt[3]{-\frac{\xi_2(\omega)}{2} - \sqrt{\frac{\xi_2^2(\omega)}{4} + \frac{\xi_1^3(\omega)}{27}}} \quad \text{when}$$

$$\frac{\xi_2^2(\omega)}{4} + \frac{\xi_1^3(\omega)}{27} > 0.$$

<u>Proposition</u>: If $x \geq y$, then

$$P(A < x, B < y, \mathcal{D}) = P\left(-x^3 - y^3 < \xi_2 < -2y^3, 3y\sqrt[3]{y^3 + \xi_2} < \xi_1 < 3x\sqrt[3]{x^3 + \xi_2}\right)$$

$$+ P\left(-2y^3 \leq \xi_2 < \infty, -3\sqrt[3]{\frac{\xi_2^2}{4}} < \xi_1 < 3x\sqrt[3]{x^3 + \xi_2}\right)$$

<u>Proof</u>: Suppose $x \geq y$. Then

$$P(A < x, B < y, \mathcal{D}) = P(A^3 < x^3, B^3 < y^3, \mathcal{D})$$

$$= P\left(-\frac{\xi_2}{2} + \sqrt{\frac{\xi_2^2}{4} + \frac{\xi_1^3}{27}} < x^3, -\frac{\xi_2}{2} - \sqrt{\frac{\xi_2^2}{4} + \frac{\xi_1^3}{27}} < y^3, \mathcal{D}\right)$$

$$= P\left(\sqrt{\frac{\xi_2^2}{4} + \frac{\xi_1^3}{27}} < x^3 + \frac{\xi_2}{2}, -\frac{\xi_2}{2} - y^3 < \sqrt{\frac{\xi_2^2}{4} + \frac{\xi_1^3}{27}}, \mathcal{D}\right)$$

$$= P\left(-\frac{\xi_2}{2} - y^3 \leq 0, \sqrt{\frac{\xi_2^2}{4} + \frac{\xi_1^3}{27}} < x^3 + \frac{\xi_2}{2}, -\frac{\xi_2}{2} - y^3 < \sqrt{\frac{\xi_2^2}{4} + \frac{\xi_1^3}{27}}, \mathcal{D}\right)$$

$$+ P\left(-\frac{\xi_2}{2} - y^3 > 0, \sqrt{\frac{\xi_2^2}{4} + \frac{\xi_1^3}{27}} < x^3 + \frac{\xi_2}{2}, -\frac{\xi_2}{2} - y^3 < \sqrt{\frac{\xi_2^2}{4} + \frac{\xi_1^3}{27}}, \mathcal{D}\right)$$

$$= P\left(-\frac{\xi_2}{2} - y^3 \leq 0, \sqrt{\frac{\xi_2^2}{4} + \frac{\xi_1^3}{27}} < x^3 + \frac{\xi_2}{2}, \mathcal{D}\right)$$

$$+ P\left(-\frac{\xi_2}{2} - y^3 > 0, \sqrt{\frac{\xi_2^2}{4} + \frac{\xi_1^3}{27}} < x^3 + \frac{\xi_2}{2}, \left(-\frac{\xi_2}{2} - y^3\right)^2 < \frac{\xi_2^2}{4} + \frac{\xi_1^3}{27}, \mathcal{D}\right)$$

$$= P\left(\xi_2 + 2y^3 \geq 0, \frac{\xi_2^2}{4} + \frac{\xi_1^3}{27} < \left(x^3 + \frac{\xi_2}{2}\right)^2, x^3 + \frac{\xi_2}{2} > 0, \mathcal{D}\right)$$

$$+ P\left(\xi_2 + 2y^3 < 0, \frac{\xi_2^2}{4} + \frac{\xi_1^3}{27} < \left(x^3 + \frac{\xi_2}{2}\right)^2, \left(\frac{\xi_2}{2} + y^3\right)^2 < \frac{\xi_2^2}{4} + \frac{\xi_1^3}{27}, x^3 + \frac{\xi_2}{2} > 0, \mathcal{D}\right)$$

$$= P\left(\xi_2 + 2y^3 \geq 0, \frac{\xi_1^3}{27} < x^6 + \xi_2 x^3, \xi_2 > -2x^3, \mathcal{D}\right)$$

$$+ P\left(\xi_2 < -2y^3, \frac{\xi_1^3}{27} < x^6 + \xi_2 x^3, y^6 + \xi_2 y^3 < \frac{\xi_1^3}{27}, \xi_2 > -2x^3, \mathcal{D}\right)$$

$$= P\left(\xi_2 \geq -2y^3, \xi_2 > -2x^3, \xi_1 < 3x\sqrt[3]{x^3 + \xi_2}, \mathcal{D}\right)$$

$$+ P\left(-2x^3 < \xi_2 < -2y^3, 3y\sqrt[3]{y^3 + \xi_2} < \xi_1 < 3x\sqrt[3]{x^3 + \xi_2}, \mathcal{D}\right)$$

Now $3y\sqrt[3]{y^3 + \xi_2} < 3x\sqrt[3]{x^3 + \xi_2}$

$$\Leftrightarrow 27y^3(y^3 + \xi_2) < 27x^3(x^3 + \xi_2)$$

$$\Leftrightarrow (y^3)^2 + \xi_2(y^3) < (x^3)^2 + \xi_2(x^3)$$

$$\Leftrightarrow (y^3)^2 + \xi_2(y^3) + \frac{\xi_2^2}{4} < (x^3)^2 + \xi_2(x^3) + \frac{\xi_2^2}{4}$$

$$\Leftrightarrow \left(y^3 + \frac{\xi_2}{2}\right)^2 < \left(x^3 + \frac{\xi_2}{2}\right)^2$$

$$\Leftrightarrow \left|y^3 + \frac{\xi_2}{2}\right| < \left|x^3 + \frac{\xi_2}{2}\right|$$

$$\Leftrightarrow \left|\xi_2 + 2y^3\right| < \left|\xi_2 + 2x^3\right|$$

$\Leftrightarrow \xi_2$ is farther from $-2x^3$ than from $-2y^3$.

Since $x \leq y$, then $-2x^3 \leq -2y^3$ so $3y\sqrt[3]{y^3 + \xi_2} < 3x\sqrt[3]{x^3 + \xi_2}$ only if

$-x^3 - y^3 < \xi_2$. Thus we continue the above equalities:

$$= P\left(\xi_2 \geq -2y^3, \xi_2 > -2x^3, \xi_1 < 3x\sqrt[3]{x^3+\xi_2}, \mathcal{D}\right)$$
$$+ P\left(-2x^3 < \xi_2 < -2y^3, -x^3-y^3 < \xi_2, 3y\sqrt[3]{y^3+\xi_2} < \xi_1 < 3x\sqrt[3]{x^3+\xi_2}, \mathcal{D}\right)$$

Now since $x \leq y$, it follows that $-2x^3 \leq -x^3 - y^3$ so above

$$= P\left(\xi_2 \geq -2y^3, \xi_2 > -2x^3, \xi_1 < 3x\sqrt[3]{x^3+\xi_2}, \mathcal{D}\right)$$
$$+ P\left(-x^3-y^3 < \xi_2 < -2y^3, 3y\sqrt[3]{y^3+\xi_2} < \xi_1 < 3x\sqrt[3]{x^3+\xi_2}, \mathcal{D}\right)$$

Now $\mathcal{D}$ is equivalent to the event $\dfrac{\xi_2^2}{4} + \dfrac{\xi_1^3}{27} > 0$, which in turn is equivalent to $\xi_1 > -3\sqrt[3]{\dfrac{\xi_2^2}{4}}$ so above

$$= P\left(\xi_2 \geq -2y^3, \xi_2 > -2x^3, -3\sqrt[3]{\dfrac{\xi_2^2}{4}} < \xi_1 < 3x\sqrt[3]{x^3+\xi_2}, \mathcal{D}\right)$$
$$+ P\left(-x^3-y^3 < \xi_2 < -2y^3, \max\left(-3\sqrt[3]{\dfrac{\xi_2^2}{4}}, 3y\sqrt[3]{y^3+\xi_2}\right) < \xi_1 < 3x\sqrt[3]{x^3+\xi_2}, \mathcal{D}\right)$$

Now we claim that $-3\sqrt[3]{\dfrac{\xi_2^2}{4}} \leq 3y\sqrt[3]{y^3+\xi_2}$. For suppose that

$-3\sqrt[3]{\dfrac{\xi_2^2}{4}} > 3y\sqrt[3]{y^3+\xi_2}$ for some $y, \xi_2$. Then

$$-27\left(\dfrac{\xi_2^2}{4}\right) > 27y^3\left(y^3+\xi_2\right) \Rightarrow 0 > 27y^3\left(y^3+\xi_2\right) + 27\left(\dfrac{\xi_2^2}{4}\right) \Rightarrow 0 > \left(y^3 + \dfrac{\xi_2}{2}\right)^2$$

which is a clear absurdity.

Thus $\max\left(-3\sqrt[3]{\dfrac{\xi_2^2}{4}}, 3y\sqrt[3]{y^3+\xi_2}\right) = 3y\sqrt[3]{y^3+\xi_2}$ so above

$$= P\left(\xi_2 \geq -2y^3, \xi_2 > -2x^3, -3\sqrt[3]{\frac{\xi_2^2}{4}} < \xi_1 < 3x\sqrt[3]{x^3+\xi_2}\right)$$

$$+ P\left(-x^3-y^3 < \xi_2 < -2y^3, 3y\sqrt[3]{y^3+\xi_2} < \xi_1 < 3x\sqrt[3]{x^3+\xi_2}\right)$$

$$= P\left(\xi_2 \geq \max(-2y^3,-2x^3), -3\sqrt[3]{\frac{\xi_2^2}{4}} < \xi_1 < 3x\sqrt[3]{x^3+\xi_2}\right)$$

$$+ P\left(-x^3-y^3 < \xi_2 < -2y^3, 3y\sqrt[3]{y^3+\xi_2} < \xi_1 < 3x\sqrt[3]{x^3+\xi_2}\right)$$

Now $x \geq y$, so $-2y^3 \geq -2x^3$, hence $\max(-2y^3,-2x^3) = -2y^3$. So above

$$= P\left(-2y^3 \leq \xi_2, -3\sqrt[3]{\frac{\xi_2^2}{4}} < \xi_1 < 3x\sqrt[3]{x^3+\xi_2}\right)$$

$$+ P\left(-x^3-y^3 < \xi_2 < -2y^3, 3y\sqrt[3]{y^3+\xi_2} < \xi_1 < 3x\sqrt[3]{x^3+\xi_2}\right)$$

and the proposition is proved.

<u>Claim</u>: $P\left(-2y^3 \leq \xi_2, -3\sqrt[3]{\frac{\xi_2^2}{4}} < \xi_1 < 3x\sqrt[3]{x^3+\xi_2}\right) = \int_{-2y^3}^{\infty} \int_{-3\sqrt[3]{\frac{u^2}{4}}}^{3x\sqrt[3]{x^3+\xi_2}} f(v,u)\,dvdu$

Proof: We need only to show that $-3\sqrt[3]{\frac{u^2}{4}} \leq 3x\sqrt[3]{x^3+u}$ for $u \in [-2y^3, \infty)$. Note that we always have the inequality

$\left(x^3+\frac{u}{2}\right)^2 \geq 0$, hence $x^6+ux^3+\frac{u^2}{4} \geq 0$, $x^3(x^3+u) \geq -\frac{u^2}{4}$, $x\sqrt[3]{x^3+u} \geq -\sqrt[3]{\frac{u^2}{4}}$,

$3x\sqrt[3]{x^3+u} \geq -3\sqrt[3]{\frac{u^2}{4}}$ as desired.

<u>Claim</u>:

$$P\left(-x^3-y^3 \leq \xi_2 < -2y^3, 3y\sqrt[3]{y^3+\xi_2} < \xi_1 < 3x\sqrt[3]{x^3+\xi_2}\right) = \int_{-x^3-y^3}^{-2y^3} \int_{3y\sqrt[3]{y^3+u}}^{3x\sqrt[3]{x^3+u}} f(v,u)\,dvdu$$

Proof: As before, we need only show that $-x^3 - y^3 \leq \xi_2 < -2y^3$ and that $-x^3 - y^3 \leq \xi_2 < -2y^3$ implies $3y\sqrt[3]{y^3+u} \leq 3x\sqrt[3]{x^3+u}$. Since $x \geq y$, then $x^3 \geq y^3$, and thus $-y^3 \geq -x^3$, implying $-x^3 - y^3 \leq -2y^3$. Next, if $-x^3 - y^3 \leq u \leq -2y^3$, then $y^3 - x^3 \leq u + 2y^3 \leq 0$, hence

$$|u+2y^3| = -(u+2y^3) = -u - 2y^3 \leq -(-x^3 - y^3) - 2y^3 = x^3 - y^3 = -x^3 - y^3 + 2x^3 \leq u + 2x^3$$

Hence $|u+2y^3| \leq u + 2x^3$

$$\Rightarrow (u+2y^3)^2 \leq (u+2x^3)^2$$

$$\Rightarrow \left(y^3 + \frac{u}{2}\right)^2 \leq \left(x^3 + \frac{u}{2}\right)^2$$

$$\Rightarrow y^6 + uy^3 + \frac{u^2}{4} \leq x^6 + ux^3 + \frac{u^2}{4}$$

$$\Rightarrow y^3(y^3 + u) \leq x^3(x^3 + u)$$

$$\Rightarrow y\sqrt[3]{y^3+u} \leq x\sqrt[3]{x^3+u}$$

$$\Rightarrow 3y\sqrt[3]{y^3+u} \leq 3x\sqrt[3]{x^3+u}, \text{ as desired.}$$

Thus we have proved the following

Proposition: If $x \geq y$, then

$$P(A<x, B<y, \mathcal{D}) = \int_{-x^3-y^3}^{-2y^3} \int_{3y\sqrt[3]{y^3+u}}^{3x\sqrt[3]{x^3+u}} f(v,u)\,dv\,du + \int_{-2y^3}^{\infty} \int_{-3\sqrt[3]{\frac{u^2}{4}}}^{3x\sqrt[3]{x^3+u}} f(v,u)\,dv\,du$$

We now use this last proposition to determine $P(A<x, B<y, \mathcal{D})$ for $x<y$. Note that $A(\omega) \geq B(\omega)$ for every $\omega \in \Omega$, hence $\{\omega; A < x, B < y, \mathcal{D}\} = \{\omega; A < x, B < x, \mathcal{D}\}$ when $x < y$. Thus if $x < y$, then $P(A < x, B < y, \mathcal{D}) = P(A < x, B < x, \mathcal{D})$. From the previous proposition, we see that $P(A < x, B < x, \mathcal{D}) = \int_{-2x^3}^{\infty} \int_{-3\sqrt[3]{\frac{u^2}{4}}}^{3x\sqrt[3]{x^3+u}} f(v,u)\,dv\,du$ which gives the next

<u>Proposition</u>: If $x < y$, then $P(A < x, B < y, \mathcal{D}) = \int_{-2x^3}^{\infty} \int_{-3\sqrt[3]{\frac{u^2}{4}}}^{3x\sqrt[3]{x^3+u}} f(v,u)\,dv\,du$

<u>Proposition</u>: If $x < y$, then $\dfrac{\partial}{\partial y} P(A < x, B < y, \mathcal{D}) = 0$

<u>Proof</u>: The previous proposition shows that if $x < y$, then $P(A < x, B < y, \mathcal{D})$ is independent of $y$.

<u>Proposition</u>: If $x > y$, then

$$\frac{\partial}{\partial y} P(A < x, B < y, \mathcal{D}) = \int_{-x^3-y^3}^{-2y^3} -\frac{6y^3 + 3u}{(y^3+u)^{\frac{2}{3}}} f\!\left(3y\sqrt[3]{y^3+u},u\right) du$$

<u>Proof</u>:

$$\frac{\partial}{\partial y} P(A < x, B < y, \mathcal{D}) = \frac{\partial}{\partial y} \int_{-x^3-y^3}^{-2y^3} \int_{3y\sqrt[3]{y^3+u}}^{3x\sqrt[3]{x^3+u}} f(v,u)\,dv\,du + \frac{\partial}{\partial y} \int_{-2y^3}^{\infty} \int_{-3\sqrt[3]{\frac{u^2}{4}}}^{3x\sqrt[3]{x^3+u}} f(v,u)\,dv\,du$$

$$= \int_{-x^3-y^3}^{-2y^3} -\frac{\partial\!\left(3y\sqrt[3]{y^3+u}\right)}{\partial y} f\!\left(3y\sqrt[3]{y^3+u},u\right) du + \frac{\partial(-2y^3)}{\partial y} \int_{3y\sqrt[3]{y^3+(-2y^3)}}^{3x\sqrt[3]{x^3+(-2y^3)}} f(v,-2y^3)\,dv$$

$$- \frac{\partial}{\partial y}(-x^3-y^3) \int_{3y\sqrt[3]{y^3+(-x^3-y^3)}}^{3x\sqrt[3]{x^3+(-x^3-y^3)}} f(v,-x^3-y^3)\,dv - \frac{\partial(-2y^3)}{\partial y} \int_{-3\sqrt[3]{\frac{(-2y^3)^2}{4}}}^{3x\sqrt[3]{x^3-2y^3}} f(v,-2y^3)\,dv$$

$$= \int_{-x^3-y^3}^{-2y^3} -\frac{6y^3+3u}{(y^3+u)^{\frac{2}{3}}} f\left(3y\sqrt[3]{y^3+u},u\right) du - 6y^2 \int_{-3y^2}^{3x\sqrt[3]{x^3-2y^3}} f(v,-2y^3) dv - 0$$

$$-(-6y^2) \int_{-3y^2}^{3x\sqrt[3]{x^3-2y^3}} f(v,-2y^3) dv$$

$$= \int_{-x^3-y^3}^{-2y^3} -\frac{6y^3+3u}{(y^3+u)^{\frac{2}{3}}} f\left(3y\sqrt[3]{y^3+u},u\right) du, \text{ as desired.}$$

<u>Proposition</u>: If $x < y$, then $\dfrac{\partial^2}{\partial x \partial y} P(A<x, B<y, \mathcal{D}) = 0$

<u>Proof</u>: $\dfrac{\partial^2}{\partial x \partial y} P(A<x, B<y, \mathcal{D}) = \dfrac{\partial}{\partial x}\left[\dfrac{\partial}{\partial y} P(A<x, B<y, \mathcal{D})\right] = \dfrac{\partial}{\partial x}[0] = 0$.

<u>Proposition</u>: If $x > y$, then

$$\dfrac{\partial^2}{\partial x \partial y} P(A<x, B<y, \mathcal{D}) = 9(x^3-y^3) f(-3xy, -x^3-y^3).$$

<u>Proof</u>:

$$\dfrac{\partial^2}{\partial x \partial y} P(A<x, B<y, \mathcal{D}) = \dfrac{\partial}{\partial x}\left\{\dfrac{\partial}{\partial y} P(A<x, B<y, \mathcal{D})\right\}$$

$$= \dfrac{\partial}{\partial x} \int_{-x^3-y^3}^{-2y^3} -\frac{6y^3+3u}{(y^3+u)^{\frac{2}{3}}} f\left(3y\sqrt[3]{y^3+u},u\right) du$$

$$= -\dfrac{\partial(-x^3-y^3)}{\partial x}\left\{-\frac{6y^3+3(-x^3-y^3)}{(y^3+(-x^3-y^3))^{\frac{2}{3}}}\right\} f\left(3y\sqrt[3]{y^3+(-x^3-y^3)},(-x^3-y^3)\right)$$

$$= -(-3x^2)\left\{-\frac{3y^3-3x^3}{x^2}\right\} f(-3xy, -x^3-y^3)$$

$$= 9(x^3-y^3) f(-3xy, -x^3-y^3)$$

Thus the conditional density of A and B with respect to $\mathcal{D}$ is given by

$$g_{A,B}(x,y|\mathcal{D}) = \begin{cases} \dfrac{9}{P(\mathcal{D})}(x^3-y^3)f(-3xy,-x^3-y^3) & \text{if } x > y \\ 0 & \text{if } x \leq y \end{cases}$$

Next, we calculate the joint density of $R_1^*$ and $R_2^*$ under the hypothesis $\mathcal{D}$.

<u>Proposition</u>: The joint density

$$h(x,y|\mathcal{D}) = \begin{cases} \dfrac{4}{P(\mathcal{D})}(y^3+9x^2y)f(y^2-3x^2,2xy^2+2x^3) & \text{if } y > 0 \\ 0 & \text{if } y \leq 0 \end{cases}$$

<u>Proof</u>: Suppose $x$ and $y$ are real numbers. Then

$$P(R_1^* < x, R_2^* < y | \mathcal{D}) = P\left(-\dfrac{A+B}{2} < x, \dfrac{(A-B)\sqrt{3}}{2} < y \middle| \mathcal{D}\right)$$

$$= P\left(A+B > -2x, A-B < \dfrac{2}{\sqrt{3}}y \middle| \mathcal{D}\right)$$

$$= P\left(B > -2x-A, A - \dfrac{2}{\sqrt{3}}y < B \middle| \mathcal{D}\right)$$

$$= P\left(\max\left(-2x-A, A - \dfrac{2}{\sqrt{3}}y\right) < B \middle| \mathcal{D}\right)$$

Now $-2x - A < A - \dfrac{2}{\sqrt{3}}y$

$$\Leftrightarrow -2x + \dfrac{2}{\sqrt{3}}y < 2A$$

$$\Leftrightarrow \dfrac{y}{\sqrt{3}} - x < A$$

So above

$$= P\left(\frac{y}{\sqrt{3}} - x < A, A - \frac{2}{\sqrt{3}}y < B < A \Big| \mathcal{D}\right) + P\left(\frac{y}{\sqrt{3}} - x \geq A, -2x - A < B < A \Big| \mathcal{D}\right)$$

$$= P\left(\frac{y}{\sqrt{3}} - x < A, A - \frac{2}{\sqrt{3}}y < B < A, y > 0 \Big| \mathcal{D}\right)$$
$$+ P\left(\frac{y}{\sqrt{3}} - x \geq A, -2x - A < B < A, -2x - A < A \Big| \mathcal{D}\right)$$

$$= P\left(\frac{y}{\sqrt{3}} - x < A, A - \frac{2}{\sqrt{3}}y < B < A, y > 0 \Big| \mathcal{D}\right) + P\left(\frac{y}{\sqrt{3}} - x \geq A > -x, -2x - A < B < A \Big| \mathcal{D}\right)$$

$$= (y > 0) P\left(\frac{y}{\sqrt{3}} - x < A, A - \frac{2}{\sqrt{3}}y < B < A \Big| \mathcal{D}\right)$$
$$+ (y > 0) P\left(\frac{y}{\sqrt{3}} - x \geq A > -x, -2x - A < B < A \Big| \mathcal{D}\right)$$

$$= (y > 0) \int_{-x}^{\frac{y}{\sqrt{3}} - x} \int_{-2x-r}^{r} g_{A,B}(r, s | \mathcal{D}) \, ds \, dr + (y > 0) \int_{\frac{y}{\sqrt{3}} - x}^{\infty} \int_{r - \frac{2y}{\sqrt{3}}}^{r} g_{A,B}(r, s | \mathcal{D}) \, ds \, dr$$

Thus if $y \neq 0$, we have

$$\frac{\partial P(R_1^* < x, R_2^* < y | \mathcal{D})}{\partial y} = (y > 0) \frac{1}{\sqrt{3}} \int_{-\frac{y}{\sqrt{3}} - x}^{\frac{y}{\sqrt{3}} - x} g_{A,B}\left(\frac{y}{\sqrt{3}} - x, s \Big| \mathcal{D}\right) ds$$
$$+ (y > 0) \left(-\frac{1}{\sqrt{3}}\right) \int_{-\frac{y}{\sqrt{3}} - x}^{\frac{y}{\sqrt{3}} - x} g_{A,B}\left(\frac{y}{\sqrt{3}} - x, s \Big| \mathcal{D}\right) ds$$
$$+ (y > 0) \int_{\frac{y}{\sqrt{3}} - x}^{\infty} -\left(-\frac{2}{\sqrt{3}}\right) g_{A,B}\left(r, r - \frac{2y}{\sqrt{3}} \Big| \mathcal{D}\right) dr$$

$$= (y > 0) \frac{2}{\sqrt{3}} \int_{\frac{y}{\sqrt{3}} - x}^{\infty} g_{A,B}\left(r, r - \frac{2y}{\sqrt{3}} \Big| \mathcal{D}\right) dr$$

This implies that if $y \neq 0$,

$$\frac{\partial P(R_1^* < x, R_2^* < y | \mathcal{D})}{\partial x \partial y} = (y > 0) \frac{2}{\sqrt{3}} g_{A,B}\left(\frac{y}{\sqrt{3}} - x, -\frac{y}{\sqrt{3}} - x \Big| \mathcal{D}\right)$$

Now

$$g_{A,B}\left(\tfrac{y}{\sqrt{3}}-x,-\tfrac{y}{\sqrt{3}}-x\Big|\mathcal{D}\right)$$

$$=\begin{cases}\dfrac{2\sqrt{3}}{P(\mathcal{D})}\left(y^3+9x^2y\right)f\left(y^2-3x^2,2xy^2+2x^3\right) & \text{if } \tfrac{y}{\sqrt{3}}-x>-\tfrac{y}{\sqrt{3}}-x \\ 0 & \text{if } \tfrac{y}{\sqrt{3}}-x\leq -\tfrac{y}{\sqrt{3}}-x\end{cases}$$

$$=\begin{cases}\dfrac{2\sqrt{3}}{P(\mathcal{D})}\left(y^3+9x^2y\right)f\left(y^2-3x^2,2xy^2+2x^3\right) & \text{if } y>0 \\ 0 & \text{if } y<0\end{cases}$$

Hence

$$h(x,y|\mathcal{D})=\begin{cases}\dfrac{4}{P(\mathcal{D})}\left(y^3+9x^2y\right)f\left(y^2-3x^2,2xy^2+2x^3\right) & \text{if } y>0 \\ 0 & \text{if } y\leq 0\end{cases}$$

and the proposition is proved.

# CALCULATION OF CONDITIONAL DENSITY $h(x,y|\mathcal{K})$

Definition:

Let $\mathcal{K} = \{(a,b) \in \mathbb{R}^2; z^3 + az + b = 0 \text{ has 3 distinct real roots}\}$

$$= \left\{(a,b) \in \mathbb{R}^2; \frac{b^2}{4} + \frac{a^3}{27} < 0\right\}.$$

Note: We assume in the following that $P((\xi_1, \xi_2) \in \mathcal{K}) > 0$.

Proposition: If $(a_0, b_0) \in \mathcal{K}$ and $t_0 \in \mathbb{R}$ is a root of $z^3 + a_0 z + b_0 = 0$, then there exists an open $(a_0, b_0)$-neighborhood $U$ and a function $R: U \to \mathbb{R}$ such that

1) $R(a_0, b_0) = t_0$;

2) $R$ has continuous first partial derivatives and is continuous in the interior of $U$;

3) For each $(a,b) \in U$, $R(a,b)$ is a root of $z^3 + az + b = 0$;

4) $R$ is the only function on $U$ with all of the above properties.

Proof: Let $F: \mathbb{R}^3 \to \mathbb{R}$ be given by $F(a,b,t) = t^3 + at + b$. Then $\frac{\partial F}{\partial t}(a,b,t) = 3t^2 + a$, $\frac{\partial F}{\partial b}(a,b,t) = 1$, and $\frac{\partial F}{\partial a}(a,b,t) = t$. Thus $F$, $\frac{\partial F}{\partial t}$, $\frac{\partial F}{\partial b}$, and $\frac{\partial F}{\partial a}$ are all continuous near the point $(a_0, b_0, t_0)$. By assumption, $t_0^3 + a_0 t_0 + b_0 = 0$, so $F(a_0, b_0, t_0) = 0$. Also, $(a_0, b_0) \in \mathcal{K}$ so the polynomial $z^3 + a_0 z + b_0 = 0$ has no multiple root. Therefore it shares no zeros with its derivative, which is $3z^2+a_0$.

But $t_0$ is a zero of $z^3+a_0z+b_0=0$, so $3t_0^2+a_0 \neq 0$, hence $\frac{\partial F}{\partial t}(a_0,b_0,t_0) \neq 0$. So by the Implicit Function Theorem, there exists a positive number $h$ and a function $s$ on the open set $B=\{(a,b)\in \mathbb{R}^2; |a-a_0|<h, |b-b_0|<h\}$ such that $s(a_0,b_0)=t_0$, $s$ has continuous first partial derivatives and is continuous on $B$, $F(a,b,s(a,b))=0$ when $(a,b)\in B$, and such that $s$ is the only such function on $B$. Take $U=B$ and $R=s$, and the proposition is proved.

Remark: $K$ is open in $\mathbb{R}^2$ and thus we may find an open neighborhood $U$ in $K$ and a function $f:U \to \mathbb{R}$ with properties 1) - 4) above if $(a_0,b_0)\in K$ and $t_0 \in \mathbb{R}$ satisfies $z^3+a_0z+b_0=0$.

Proposition: Let $R_1, R_2, R_3: K \to \mathbb{R}$ be given by the conditions $R_1(a,b) > R_2(a,b) > R_3(a,b)$ for $(a,b)\in K$ and $[R_i(a,b)]^3 + a[R_i(a,b)] + b = 0$ for $(a,b)\in K$ and $i=1,2,3$. Then $R_1$, $R_2$, and $R_3$ are continuously differentiable on $K$.

Proof: Fix $(a_0,b_0)\in K$. There exists an open neighborhood $V \subseteq K$ such that $(a_0,b_0)\in V$ and unique functions $S_1$, $S_2$, and $S_3$ on $V$ such that for $i=1,2,3$, that $S_i(a_0,b_0)=R_i(a_0,b_0)$ and $S_i$ is continuously differentiable on $V$, and such that $S_i(a,b)$ is a root of $z^3 + az + b = 0$ for $(a,b) \in V$, and such that $S_1(a,b)$, $S_2(a,b)$, and $S_3(a,b)$ are distinct for every $(a,b) \in V$. Since they are continuous on $V$ and $S_1(a_0,b_0) > S_2(a_0,b_0) > S_3(a_0,b_0)$, we see that there exists an open $(a_0,b_0)$-nbhd $W \subseteq V$ such that $S_1(a,b) > S_2(a,b) > S_3(a,b)$ for $(a,b) \in W$. Thus $R_1 \equiv S_1$, $R_2 \equiv S_2$, $R_3 \equiv S_3$ on $W$, so $R_1$, $R_2$, and $R_3$ are continuously differentiable at $(a_0,b_0)$, hence on $K$.

<u>Definition</u>: Let $g: K \to \mathbb{R}^2$ be given by $g(u,v) = (R_1(u,v), R_2(u,v))$

<u>Proposition</u>: $Rng(g) = \{(x,y); x > 0, -\frac{1}{2}x < y < x\}$.

<u>Proof</u>: Choose $(x,y) \in Rng(g)$. Let $(u,v) \in K$ be such that $g(u,v) = (x,y)$, i.e. $R_1(u,v) = x$ and $R_2(u,v) = y$. Suppose that $x \leq 0$. Then $R_1(u,v) \leq 0$, so $R_3(u,v) < R_2(u,v) < R_1(u,v) \leq 0$, hence $0 = R_3(u,v) + R_2(u,v) + R_1(u,v) < 0$, which is absurd. Hence $x > 0$.

Next, $y = R_2(u,v) < R_1(u,v) = x$ so $y < x$. Next, suppose that $-\frac{1}{2}x \geq y$. Then $R_2(u,v) = y \leq -\frac{1}{2}x \leq -x - y = -R_1(u,v) - R_2(u,v) = R_3(u,v) < R_2(u,v)$, which is absurd. Hence $-\frac{1}{2}x < y < x$.

Conversely, suppose that $(x,y)$ is such that $x > 0$ and $-\frac{1}{2}x < y < x$. Note that $z^3 - (x^2 + xy + y^2)z + xy(x+y) \equiv (z-x)(z-y)(z+x+y)$ has zeros $x$, $y$, $-x-y$. Since $-\frac{1}{2}x < y < x$, we see that $x > y > -x-y$. Since the roots are distinct, it follows that $(-x^2 - xy - y^2, xy(x+y)) \in K$.

Hence $R_1(-x^2 - xy - y^2, xy(x+y)) = x$ and $R_2(-x^2 - xy - y^2, xy(x+y)) = y$, thus $g(-x^2 - xy - y^2, xy(x+y)) = (x, y)$, hence $(x, y) \in \text{Rng}(g)$, and the proposition is proved.

<u>Proposition</u>: $g: K \to \mathbb{R}^2$ is continuously differentiable on $K$.

<u>Proof</u>: Follows immediately from the fact that $R_1$ and $R_2$ are continuously differentiable on $K$.

<u>Proposition</u>: $g: K \to \mathbb{R}^2$ is 1 - 1 on $K$.

<u>Proof</u>: Suppose $g(u_0, v_0) = g(u_1, v_1)$ for some $(u_0, v_0), (u_1, v_1) \in K$. Then $R_1(u_0, v_0) = R_1(u_1, v_1)$ and $R_2(u_0, v_0) = R_2(u_1, v_1)$. Next,

$R_3(u_0, v_0) = -R_1(u_0, v_0) - R_2(u_0, v_0) = -R_1(u_1, v_1) - R_2(u_1, v_1) = R_3(u_1, v_1)$. Next,

$u_0 = R_1(u_0, v_0) R_2(u_0, v_0) + R_1(u_0, v_0) R_3(u_0, v_0) + R_2(u_0, v_0) R_3(u_0, v_0)$
$= R_1(u_1, v_1) R_2(u_1, v_1) + R_1(u_1, v_1) R_3(u_1, v_1) + R_2(u_1, v_1) R_3(u_1, v_1) = u_1$

Also, $v_0 = -R_1(u_0, v_0) R_2(u_0, v_0) R_3(u_0, v_0) = -R_1(u_1, v_1) R_2(u_1, v_1) R_3(u_1, v_1) = v_1$

Hence $(u_0, v_0) = (u_1, v_1)$, and therefore $g$ is 1 - 1 on $K$.

<u>Proposition</u>: $g$ is invertible on $g(K)$, and $g^{-1}: g(K) \to K$ is given by $g^{-1}(x, y) = (-x^2 - xy - y^2, xy(x+y))$.

Proof: Suppose $(x,y) \in g(K)$. Then there exists a unique $(u,v) \in K$ such that $R_1(u,v) = x$ and $R_2(u,v) = y$. Then $R_3(u,v) = -R_1(u,v) - R_2(u,v) = -x - y$. Hence

$$u = R_1(u,v)R_2(u,v) + R_1(u,v)R_3(u,v) + R_2(u,v)R_3(u,v)$$
$$= (x)(y) + (x)(-x-y) + (y)(-x-y) = -x^2 - xy - y^2$$

and

$$v = -R_1(u,v)R_2(u,v)R_3(u,v) = -(x)(y)(-x-y) = xy(x+y).$$

Remark: $g^{-1}$ is continuously differentiable on $g(K)$.

Proposition: The Jacobian of the transformation $g^{-1}: g(K) \to K$ is the map $J: g(K) \to \mathbb{R}$ given by $J(x,y) = -(x-y)(2x+y)(x+2y)$.

Proof: $J(x,y) = \begin{vmatrix} \dfrac{\partial g_1^{-1}}{\partial x} & \dfrac{\partial g_1^{-1}}{\partial y} \\ \dfrac{\partial g_2^{-1}}{\partial x} & \dfrac{\partial g_2^{-1}}{\partial y} \end{vmatrix} = \begin{vmatrix} -2x-y & -x-2y \\ 2xy+y^2 & x^2+2xy \end{vmatrix}$

$$= (-2x-y)(x^2+2xy) - (-x-2y)(2xy+y^2)$$
$$= (-2x-y)(x)(x+2y) - (-x-2y)(y)(2x+y) = -(x-y)(2x+y)(x+2y)$$

Proposition: $0 \notin \mathrm{Rng}(J)$

Proof: Fix $(x,y) \in g(K)$. Then $x > 0$ and $-\frac{1}{2}x < y < x$. $x - y > 0$ since $x > y$. $2x + y > 0$ because $2x + y > 2x + (-\frac{1}{2}x) = \frac{3}{2}x > 0$. $x + 2y > 0$ because $x + 2y > x + 2(-\frac{1}{2}x) = 0$. Hence $J(x,y) < 0$, thus $0 \notin \mathrm{Rng}(J)$.∎

Proposition: If the joint density $f$ of $(\xi_1, \xi_2)$ is continuous on $K$ and if $S$ is any closed, bounded subset of $g(K)$ which has area, then

$$P\big((R_1(\xi_1,\xi_2), R_2(\xi_1,\xi_2)) \in S\big)$$
$$= \iint_S (x-y)(2x+y)(x+2y) f(-x^2-xy-y^2, xy(x+y)) dxdy$$

Proof: Let $f$ and $S$ be as above. Then note that

$$\{\omega \in \Omega; (R_1(\xi_1(\omega),\xi_2(\omega)), R_2(\xi_1(\omega),\xi_2(\omega))) \in S\}$$
$$= \{\omega \in \Omega; (\xi_1(\omega),\xi_2(\omega)) \in K, g(\xi_1(\omega),\xi_2(\omega)) \in S\}$$
$$= \{\omega \in \Omega; (\xi_1(\omega),\xi_2(\omega)) \in K, (\xi_1(\omega),\xi_2(\omega)) \in g^{-1}S\}$$
$$= \{\omega \in \Omega; (\xi_1(\omega),\xi_2(\omega)) \in g^{-1}S\}$$

So

$$P\big((R_1(\xi_1,\xi_2), R_2(\xi_1,\xi_2)) \in S\big) = \iint_{g^{-1}(S)} f(u,v) dudv$$
$$= \iint_S f(g_1^{-1}(x,y), g_2^{-1}(x,y)) |J(x,y)| dxdy$$

where $g^{-1}(x,y) = (g_1^{-1}(x,y), g_2^{-1}(x,y))$.

Substitution in the above integral yields

$$\iint_S f(-x^2-xy-y^2, xy(x+y))(x-y)(2x+y)(x+2y) dxdy,$$

which is the desired result.

Remark: If $f$ is continuous and $S$ is closed and bounded but does not meet $g(K)$, then $P\big((R_1(\xi_1,\xi_2), R_2(\xi_1,\xi_2)) \in S\big) = 0$. This is because

$$\{\omega \in \Omega; (R_1(\xi_1(\omega), \xi_2(\omega)), R_2(\xi_1(\omega), \xi_2(\omega))) \in S\}$$
$$= \{\omega \in \Omega; g(\xi_1(\omega), \xi_2(\omega)) \in S\}$$
$$= \{\omega \in \Omega; g(\xi_1(\omega), \xi_2(\omega)) \in S \cap g(K)\}$$
$$= \{\omega \in \Omega; g(\xi_1(\omega), \xi_2(\omega)) \in \varnothing\} = \varnothing$$

<u>Theorem</u>: If the joint density $f$ of $(\xi_1, \xi_2)$ is continuous on $K$ and if $S$ is any measurable subset of $\mathbb{R}^2$, then

$$P\big((R_1(\xi_1, \xi_2), R_2(\xi_1, \xi_2)) \in S\big)$$
$$= \iint_S \chi_{g(K)}(x, y)(x-y)(2x+y)(x+2y) f\big(-x^2 - xy - y^2, xy(x+y)\big) dx dy$$

<u>Theorem</u>: If the joint density $f$ of $(\xi_1, \xi_2)$ is continuous on $K$, then the joint conditional density $h(x, y | K)$ for $(R_1(\xi_1, \xi_2), R_2(\xi_1, \xi_2))$ is given by

$$h(x, y | \mathcal{K}) = \begin{cases} \dfrac{(x-y)(2x+y)(x+2y)}{P(\mathcal{K})} f\big(-x^2 - xy - y^2, xy(x+y)\big) & \text{if } (x, y) \in g(K) \\ 0 & \text{if } (x, y) \notin g(K) \end{cases}$$

except possibly on a set of measure zero.

SUMMARY

We have found that

$$h(x,y|\mathcal{D}) = \frac{4}{P(\mathcal{D})}(y^3 + 9x^2y)f(y^2 - 3x^2, 2xy^2 + 2x^3) \quad (y > 0)$$

$$h(x,y|\mathcal{K}) = \frac{(x-y)(2x+y)(x+2y)}{P(\mathcal{K})}f(-x^2 - xy - y^2, x^2y + xy^2) \quad (x > 0, -\tfrac{1}{2}x < y < x)$$

# VITA

Kerry Michael Soileau was born in New Orleans, Louisiana on June 8, 1956. He was graduated from Florida Technological University with a B.S. in mathematics on June 11, 1976. He attended Duke University during the 1976-77 academic year and was later employed by Amoco Production Company as an exploration technologist, and by the National Aeronautics and Space Administration as an Aerospace Summer Intern. At present he is a candidate for the degree of Master of Science in the Department of mathematics at Louisiana State University.